\documentclass[12pt]{article}
\usepackage{amssymb}
\usepackage{enumerate}
\usepackage{graphicx}
\usepackage{amsthm}
\usepackage{color}
\def\C{\mathbb{C}}

\def\R{\mathbb{R}}

\def\dst{\displaystyle}

\newtheorem{theorem}{Theorem}
\newtheorem{corollary}{Corollary}

\newtheorem{lemma}{Lemma}

\title{Injectivity of non-singular planar maps with  disconnecting curves in the eigenvalues space}

\author{M. Sabatini \footnote{Universit\`a di Trento, Dipartimento di Matematica, Via Sommarive 14, I-38121, Povo, Trento - Italy. Email: marco.sabatini@unitn.it. This work has been partially supported by GNAMPA. MSC: 14R15} }
\begin{document}
\maketitle
\begin{abstract}  Fessler and Gutierrez \cite{Fe,Gu} proved that if a non-singular planar map has Jacobian matrix without   eigenvalues in $(0,+\infty)$, then it is injective. We prove that the same holds replacing $(0,+\infty)$ with any unbounded curve disconnecting the upper (lower) complex half-plane. Additionally we prove that a Jacobian map $(P,Q)$ is injective if $P_x + Q_y$ is not a surjective function.

 \end{abstract}

\section{Introduction}

Let us consider a map $F = (P,Q) \in C^1(\R^2,\R^2)$. Let
$$
J_F (x,y) = \left( \matrix{P_x(x,y) & P_y(x,y) \cr Q_x(x,y) & Q_y(x,y)}  \right)
$$ 
be the jacobian matrix of $F$ at $(x,y)$. We denote by $T(x,y)$ the trace of $J_F(x,y)$, i. e. the divergence of the vector field $F(x,y)$,  by $D(x,y)$ its determinant and by $\Delta(x,y) = T(x,y)^2 - 4 D(x,y)$ the discriminant of the eigenvalues equation. 
We say that $F(x,y)$ is a \textit{non-singular map} if $D(x,y) \neq 0$ on all of $\R^2$, and that it  is a \textit{Jacobian map} if $D(x,y) $ is a non-zero constant on all of $\R^2$. We denote by $\Sigma_F(x,y)$ the spectrum of $J_F(x,y)$, i. e. the set of its eigenvalues. We set
$\Sigma_F = \cup\, \{ \Sigma_F(x,y) : (x,y) \in \R^2 \}$. 

The implicit function theorem gives the  injectivity of a map in a neighbourhood of a point   $(x^*,y^*) \in \R^2$  such that  $D(x^*,y^*) \neq 0$. On the other hand, even if $D(x,y) \neq 0$ on all of $\R^2$ the map can be non-injective, as  the exponential map $(e^y \cos x,e^y \sin x)$. The search for additional conditions ensuring the global   injectivity of a locally invertible map is a classical problem. A fundamental result is Hadamard global inverse function theorem, which gives the global invertibility of a proper non-singular map $F\in C^1(\R^n,\R^n)$. 
 In this field some old  problems still resist the attempts to find  a solution. The celebrated Jacobian Conjecture is  concerned with polynomial maps $F: \C^n \to \C^n$ \cite{Ke}. According to such a conjecture, a polynomial map with non-zero constant jacobian determinant is invertible, with polynomial inverse. Such a statement and its variants were studied  in different settings, even replacing $\C^n$ with $\R^n$ or other fields, and several partial results were proved, but it is not yet proved or disproved even for $n=2$ \cite{BCW,vdE}. Another famous problem, known as the Global Asymptotic Stability Jacobian conjecture \cite{MY} was proved in the planar case to be equivalent to a global injectivity one.  Such a conjecture was proved to be true in dimension 2  \cite{Fe,Gl,Gu}, false in higher dimensions \cite{CGEHL}.
 
  In \cite{Fe,Gu} the injectivity of a map with $D(x,y) > 0$ was proved under the additional assumption that for $(x,y) \not \in K$, $K$ compact, the eigenvalues of $J_F(x,y)$ do not belong to $(0,+\infty)$.  Such a result was extended in \cite{CGL,GS,GV,Ra}. Other results proving injectivity with different additional conditions were obtained in \cite{BGL,BdSF,BO,Ch}.
 
 In this paper we propose an approach based on the eigenvalues continuity. We prove that the  injectivity can be proved by replacing the half-line $(0,+\infty)$ by any unbounded curve $\delta$ in the complex plane, provided $\delta$ disconnects the upper (lower) half-plane. This allows to prove the injectivity as a consequence of some suitable inequalities. Moreover we prove that a Jacobian map is injective if there exists $z \in (-\infty) \cup S^1 \cup (0,+\infty)$ such that $z$ is not an eigenvalue of $J_F(x,y)$, for $(x,y) \in \R^2$.
 As a consequence, if  the function $T(x,y)$ is not surjective, then $F(x,y)$ is injective. We do not require $F(x,y)$ to be polynomial.
  
 \section{Results}

We report next  Lemma without proof, since it is a standard statement in finite dimensional spectral theory. We denote by $Re(\lambda)$,  $Im(\lambda)$, resp. the real and imaginary part of the complex number $\lambda$.

\begin{lemma}  \label{lemma} 
Let $F\in C^1(\R^2,\R^2)$. Then there exist functions  $\lambda_1, \lambda_2 \in C^0(\R^2,\R^2)$, such that for all $(x,y) \in \R^2$   $\lambda_1(x,y)$ and $\lambda_2(x,y)$ are the eigenvalues of  $J_F(x,y)$. Such functions can be taken such that $Re(\lambda_1(x,y)) \geq 0$ and \\
$Re(\lambda_2(x,y)) \leq 0$.
\end{lemma}

Either such eigenvalues are real or complex conjugate. This implies that the set $\lambda_1(\R^2) \cup \lambda_2(\R^2) $ is symmetric with respect to the $x$ axis.   For the reader's convenience we report the main theorem proved in \cite{Fe,Gu}, that will be applied in the following.

\begin{theorem}  \label{fegu} (Fessler - Gutierrez)
Let $F\in C^1(\R^2,\R^2)$ with $D(x,y) > 0$ for all $(x,y) \in \R^2$. Assume there exists a compact set $K \subset \R^2$ such that for all $(x,y) \not\in K$ the eigenvalues of $J_F(x,y)$ are not in $(0,+\infty)$. Then  $F$ is injective \cite{Fe,Gu}.  
\end{theorem}

When convenient, in the following we sometimes identify $\C$ with the real plane $\R^2$. Let us set $\C^+ = \{ u +i v: v \geq  0 \}$.  We say that an unbounded  curve $\delta \in C^0([0,+\infty),\C^+)$ disconnects  the half-plane $\C^+$  if $\delta(0) = 0+0\cdot i \equiv (0,0)$, $\delta$ has no other points on the real axis and there exist two   connected subsets $A$, $B$, such that $\C^+ = A \cup B$, $\partial A = \delta = \partial B$. We do not require $A$ and $B$ to be disjoint. Such a definition implies that the open real half-axes  $(0,+\infty) $ and $(-\infty,0)  $ are not both contained  in $A$ or in $B$. 

We write $K^c$ for the set-theoretical complement of a set $K$.

\begin{theorem}  \label{eigencurve} 
Let $F\in C^1(\R^2,\R^2)$ with $D(x,y) > 0$ for all $(x,y) \in \R^2$. Assume there exists a compact set $K \subset \R^2$ and a  curve $\delta \in C^0([0,+\infty),\C^+)$ disconnecting $\C^+$ and such that for all $(x,y) \not \in K$ the  eigenvalues of $J_F(x,y)$ are not on $\delta$.
Then  $F$ is injective.  
\end{theorem}
{\it Proof.}  
We prove that the hypotheses of theorem \ref{fegu} either hold for $F(x,y)$ or for $-F(x,y)$. 
By absurd, assume that  neither $F(x,y)$ nor  $-F(x,y)$ satisfy them.
Hence for every compact $K \subset \R^2$ there exist both $(x_K^+,y_K^+) \not \in K$,  $(x_K^-,y_K^-)  \not \in K$ and two eigenvalues $\lambda(x_K^+,y_K^+) \in (0,+\infty)$ and $\lambda(x_K^-,y_K^-) \in (-\infty,0)$. This implies that also
 $\dst{  \frac {D(x_K^+,y_K^+)}{\lambda(x_K^+,y_K^+)} \in (0,+\infty)  }$ and   $\dst{  \frac {D(x_K^-,y_K^-)}{\lambda(x_K^-,y_K^-)} \in (-\infty,0)  }$ are eigenvalues. Hence one has $\lambda_i(x_K^+,y_K^+) \in (0,+\infty)$ and $\lambda_i(x_K^-,y_K^-) \in (-\infty,0)$, $i = 1,2$.
By the compactness of $K$ there exists a curve $\gamma \in C^0([-1,1], \R^2)$ with no points in $K$ and  connecting  $(x_K^+,y_K^+)$ to  $(x_K^-,y_K^-)$: $\gamma(-1) = (x_K^-,y_K^-)$, $\gamma(1) = (x_K^+,y_K^+)$. 

Let us consider the eigenvalue function $\lambda_1(x,y)$ with $Re(\lambda_1(x,y)) \geq 0$, as in Lemma \ref{lemma}. Let us consider the curve $\lambda_1(\gamma(t))$. One has $\lambda_1(\gamma(-1)) \in (-\infty,0)$, 
$\lambda_1(\gamma(1)) \in (0,+\infty)$. By hypothesis both points are not on $\delta $, and since $\delta(0) = (0,0)$ one of them  belongs to the set $A$, the other one to $B$. The curve $\lambda_1(\gamma(t))$ connects them, hence it has to cross the common boundary of $A$ and $B$, which is $\delta$. This contradicts the hypothesis that no eigenvalues are on $\delta$.

\hfill  $\clubsuit$  \\

Curves disconnecting $\C^+$ may be very complex. In order to deduce simple conditions for injectivity we consider a simple class of separating curves $\delta$, i. e. the graphs of the functions $v =  a u^b$, $a, b > 0$, $ u \geq 0$, or $v =  a (-u)^b$, $a, b > 0 \geq u$.  

\begin{corollary}  \label{corolla1} 
Let $F\in C^1(\R^2,\R^2)$ with $D(x,y) > 0$ for all $(x,y) \in \R^2$. Assume there exists  a compact set $K \subset \R^2$ and $a, b \in \R$, $a, b > 0$, such that for all $(x,y) \not \in K$ one of the following condition holds:
\begin{itemize}
\item[i)] if $T \geq 0$ then $\sqrt{| \Delta |} \neq a\, T^{b} $;
\item[ii)] if $T \leq 0$ then $\sqrt{| \Delta |} \neq a\, (-T)^{b} $;
\end{itemize}
then $F(x,y)$ is injective.
\end{corollary}
{\it Proof.} We prove only $i)$, the statement $ii)$ can be proved similarly. 

Since $a, b >0$, the curve $\delta$ of equation $v = a u^b$, $u \geq 0$, starts at the origin and separates $\C^+$. Such a curve  has no points on the $x$ axis, except the origin, hence if an eigenvalue $\lambda$ belongs to $\delta$ its imaginary part is not zero. This implies that $ \Delta  =  T^2 - 4D < 0$, hence   $| \Delta | =   4 D - T^2 $. 
An eigenvalue $\lambda$ belongs to $\delta$ if and only if $Re(\lambda) > 0$ and 
\begin{equation} \label{DTab}
Im(\lambda) = a\, Re(\lambda) ^b  \quad \iff \quad  \sqrt{| \Delta |} =  \sqrt{ 4 D - T^2 } = a T^b. 
\end{equation}

By hypothesis this does not occur, hence the thesis.

\hfill  $\clubsuit$  \\

We emphasize that corollary \ref{corolla1} contains two independent statements. Statement i) is not concerned with points where $T(x,y) < 0$; statement ii) is not concerned with points where $T(x,y) > 0$. 

In next corollary we prove the  injectivity under a suitable assumption on the ratio $\dst{\frac {T^2}D}$. 

\begin{corollary}  \label{corolla2} Let $F\in C^1(\R^2,\R^2)$ with $D(x,y) > 0$ for all $(x,y) \in \R^2$. Assume there exists  a compact set $K \subset \R^2$ and $c \in [0,4]$ such that for all $(x,y) \not \in K$ one has:
$$
\frac {T^2}D \neq c,   
$$
then $F(x,y)$ is injective.
\end{corollary}
{\it Proof.} If $c = 0$, then $T(x,y)$ does not vanish, has constant sign and one can apply the theorems about the Global Asymptotic Stability Jacobian Conjecture \cite{Fe,Gl,Gu}.

 If $c \in (0,4)$, we take as $\delta$ the line  of equation $\sqrt{4 - c}\, u - \sqrt{c}\, v = 0$. Assume by absurd that an eigenvalue $\lambda$ belongs to $\delta$. Then one has $\Delta < 0$ and 
$$
0 = \sqrt{4 - c} \, T - \sqrt{ c}\, \sqrt{4 D - T^2}  ,
$$
$$
 (4 - c) T^2 = 4 c D - c T^2   ,
$$
which implies $ T^2 =  c D $, contradiction.

If $c = 4$, then either $\Delta < 0$ on all of $K^c$, or $\Delta > 0$ on all of $K^c$. In the former case the eigenvalues are not real, hence theorem \ref{fegu} applies. In the latter they are real and one can take the imaginary axis as $\delta$.

\hfill  $\clubsuit$  \\

The situation is much simpler when dealing with real Jacobian maps. We can always reduce to the case $D(x,y) \equiv 1$, by possibly multiplying one component by a suitable non-zero constant.  If $D(x,y) \equiv 1$, the  eigenvalues are contained in the set
$$
\mathbb{G} =   (-\infty,0) \cup S^1 \cup (0,+\infty) ,
$$
where $S^1$ is the unit circle in $\C$. Such eigenvalues appear in couples $\dst{\lambda, \frac 1\lambda}$, if real, or $u \pm i v$, if non-real. The set $ \Sigma_F(x,y)$ is symmetric w. resp. to the real axis, i. e. it coincides with its conjugate 
$\overline{ \Sigma_F(x,y)  }$. Disconnecting $\mathbb{G}$ requires at most a couple of points. This is used in next statements in order to prove injectivity.
For the reader's convenience we report the main theorem proved in \cite{Ra}, that will be applied in the following.

 \begin{theorem}  \label{ra} (Rabanal) 
Let $F\in C^1(\R^2,\R^2)$. If there exists $\varepsilon > 0$ such that $J_F(x,y)$ has no eigenvalues in $[0,\varepsilon)$, then  $F$ is injective.  
\end{theorem}

\begin{theorem}  \label{teorema2}  Let $F\in C^1(\R^2,\R^2)$ with $D(x,y) \equiv 1$ on all of $ \R^2$. If one of the following conditions holds, then $F$ is injective:
\begin{itemize}
\item[i)]  there exists  a compact set $K \subset \R^2$ and $z \in S^1$  such that for all $(x,y) \not \in K$ one has $z \not \in \Sigma_F(x,y)$.
\item[ii)] there exists  $z \in \R \setminus \{ -1,0,1 \}  $  such that for all $(x,y)  \in \R^2$ one has $z \not \in \Sigma_F(x,y)$.

\end{itemize}
\end{theorem}
{\it Proof.} 

i)
Let $\Omega_K$ be a closed disk large enough to have $K \subset \Omega_K$. For all $(x,y) \not \in \Omega_K$ one has $z \not \in \Sigma_F(x,y)$.  The continuous maps $\lambda_i$, $i = 1,2$, map the connected set $\Omega_K^c$ into connected subsets $\lambda_i(\Omega_K^c)$, $i = 1,2$, of $\mathbb{G}$. 
We consider three cases.

i.1) If $z = u + iv \in S^1 $, $z \neq \pm 1$, is not an eigenvalue, then also $\overline{z} = u - iv  \in S^1$  is not an eigenvalue.
The couple $u \pm i v$ disconnects $\mathbb{G}$. One has 
$$
\mathbb{G} \setminus \{ u-iv, u+iv \} =  \mathbb{G}_- \cup \mathbb{G}_+ ,
$$
where $\mathbb{G}_-$ and $ \mathbb{G}_+$ are connected and $(-\infty,0) \subset \mathbb{G}_-$, $(0, + \infty) \subset \mathbb{G}_+$. 
If $\lambda_1(\Omega_K^c) \subset \mathbb{G}_-$, then also  $\lambda_2(\Omega_K^c) \subset \mathbb{G}_-$, hence $J_F(x,y)$ has  no eigenvalues in $(0,+\infty)$ for $(x,y) \not \in \Omega_K$, thus proving the injectivity of $F$. 
Similarly, if 
$\lambda_1(\Omega_K^c) \subset \mathbb{G}_+$, then also  $\lambda_2(\Omega_K^c) \subset \mathbb{G}_+$ and $J_F(x,y)$ has  no eigenvalues in $(-\infty,0)$ for $(x,y) \not \in \Omega_K$, thus proving the injectivity of $-F$, hence that one of $F$. 

i.2) $\dst{  z = \frac 1z =  1 }  $. Then the number 1 disconnects  $\mathbb{G} $ and one can write
$$
\mathbb{G} \setminus \{ 1 \} = \mathbb{G}_- \cup (0,1) \cup (1,+\infty) ,
$$
where we have set $\mathbb{G}_- =  \Big(  (-\infty,0) \cup S^1 \Big)  \setminus \{ 1 \} $.
If for some $(x,y) \in \Omega_K^c$ the matrix $J_F(x,y)$ has a positive eigenvalue, then both eigenvalues are positive and by the connectedness of $\lambda_i(\Omega_K^c)$, $i=1,2$, one has
$$
\lambda_1(\Omega_K^c) \cup \lambda_2(\Omega_K^c) \subset (0,1) \cup (1,+\infty).
$$
As a consequence $(-\infty,0)$ contains no eigenvalues, so that $-F$ is injective. 

On the other hand, if for some $(x,y) \in \Omega_K^c$ the matrix $J_F(x,y)$ has an eigenvalue in $\mathbb{G}_- $, than both eigenvalues are in $\mathbb{G}_- $ and by connectedness 
$$
\lambda_1(\Omega_K^c) \cup \lambda_2(\Omega_K^c) \subset \mathbb{G}_-  \ .
$$
Hence there are no eigenvalues in $(0,+\infty)$ and theorem \ref{fegu} gives the  injectivity of $F$.

i.3) $\dst{  z = \frac 1z = - 1 }  $. Similar to i.2).

ii) If $z \in \R \setminus \{ -1,0,1 \}  $ is not an eigenvalue, then the numbers $z$ and $\dst{ \frac 1z}$ have the same sign and disconnect the set $\mathbb{G}$ into three connected sets. For instance, if $z \in (0,1)$ we can write
$$
\mathbb{G} = (0,z) \cup \left( \frac 1z , +\infty \right) \cup \left(  \left(z,\frac 1z \right) \cup S^1  \right).
$$
Similarly, exchanging $z$ and $\dst{\frac 1 z}$, if $z \in (1,+\infty)$. As in case i.2) at least one component is free of eigenvalues for $(x,y) \not \in \Omega_K$. If $(0,z)$ does not contain eigenvalues, then one can apply theorem  \ref{ra} with $\varepsilon = z$ in order to get injectivity. Similarly if  $\dst{   \left( \frac 1z , +\infty \right)  }$ does not contain eigenvalues; such a case is equivalent to $(0,z)$  not containing eigenvalues.  If $\dst{  \left(z,\frac 1z \right) \cup S^1    }$ does not contain eigenvalues, then $(-\infty,0)$ is free of eigenvalues and applying either theorem \ref{fegu} or theorem \ref{ra} to $-F$ one proves the injectivity of $F$.

\hfill  $\clubsuit$  \\

We can deduce  a simple corollary from theorem \ref{teorema2}.

\begin{corollary}  \label{corteorema2} Let $F\in C^1(\R^2,\R^2)$ with $D(x,y) \equiv 1$ on all of $ \R^2$. If there exists $z \in (-\infty) \cup S^1 \cup (0,+\infty)$ which is not an eigenvalue of $J_F(x,y)$, for any $(x,y) \in \R^2$,
then $F(x,y)$ is injective.
\end{corollary}
{\it Proof.} 
Under the above hypothesis either i) or ii) of theorem \ref{teorema2} hold on all of $\R^2$. 
\hfill  $\clubsuit$  \\

The condition on the eigenvalues can be deduced from suitable conditions on $T(x,y)$.

\begin{corollary}  \label{trasur} Let $F\in C^1(\R^2,\R^2)$ with $D(x,y) \equiv 1$ on all of $ \R^2$. If one of the following conditions holds, then $F$ is injective:
\begin{itemize}
\item[i)]  there exists  a compact set $K \subset \R^2$ and $h \in [-2,2]$ such that for all $(x,y) \not \in K$ one has $T(x,y) \neq h$.
\item[ii)] there exists $h \in (-\infty,-2) \cup (2,+\infty)$ such that for all $(x,y) \in \R^2$, $T(x,y) \neq h$.
\end{itemize}
\end{corollary}
{\it Proof.} 
i) 
$\lambda \neq 0$ is an eigenvalue if and only if for some $(x,y)$ one has
$$
\lambda^2 - T(x,y) \lambda + 1 = 0,
$$
hence 
$$
  \lambda = \frac {T(x,y) \pm \sqrt{T(x,y)^2 - 4}}2 
  $$

If there exists $h \in \R$ such that   $T(x,y) \neq h$, then 
$\dst{  z_{1,2} = \frac {h \pm \sqrt{h^2 - 4}}2 }$
are not   eigenvalues of $J_F$. In case i) one has  $z_{1,2} \in S^1$ and point i) of theorem \ref{teorema2} applies. 

 In case ii) one has $z_{1,2}  \in  \mathbb{G} \setminus S^1 $ and point ii) of theorem \ref{teorema2} applies. 

\hfill  $\clubsuit$  \\

As a consequence we have the following corollary.

\begin{corollary}  \label{surtrace} Let $F\in C^1(\R^2,\R^2)$ with $D(x,y) \equiv 1$ on all of $ \R^2$. If  $T(x,y)$ is not surjective, then $F(x,y)$ is injective.
\end{corollary}
{\it Proof.} 
If $T(x,y)$ is not surjective, then either i) or ii) of corollary \ref{trasur} hold on all of $\R^2$. 
\hfill  $\clubsuit$  \\

The  hypotheses of corollary \ref{surtrace} do not apply to even-degree  polynomial maps. In fact, if $F(x,y)$ is an even-degree  polynomial map, then $T(x,y)$ is an odd-degree polynomial, hence it is surjective. On the other hand, odd-degree Jacobian maps with non-surjective $T(x,y)$ do exist. An example of polynomial Jacobian map with non-surjective $T(x,y)$ is given by $\dst{  F(x,y) = \big( x + y^3, y - x^3 - 3 x^2y^3 -3 x y^6 - y^9 \big)  }$. In this case one has $T(x,y) = 2 - 9 y^2 (x + y^3)^2$ which does not assume values greater than 2.

\end{document}